%% file: SL2Z_complet.tex
\input JOSEPH

%






\centerline{$\displaystyle SL(2, {\Bbb Z})$%
, les tresses \`a trois brins et le tore modulaire.}
\vskip5mm
\centerline{{\it par} Emil Artin {\it et\/} Jean Cerf}
%
\vskip3mm
\centerline{\eightsl R\'esum\'e}
{\eightpoint\sl
L'action de 
$SL(2, {\bf Z})$
 sur le tore entier et son quotient par sym\'etrie centrale
ram\`ene \`a la pr\'esentation d'Artin
des tesses \`a trois brins,
une pr\'esentation avec les g\'en\'erateurs paraboliques
$\pmatrix{1&\!\! -1\cr 0& 1\cr}$
et
$\pmatrix{1& 0\cr 1& 1\cr}$
 et explicite l'action du groupe d\'eriv\'e sur
le demi-plan de Poincar\'e.
}
\vskip3mm
\centerline{\eightsl Abstract}
{\eightpoint\sl
The action of 
$SL(2, {\bf Z})$
on the integer torus and its quotient by central symmetry
and Artin's presentation of three strings braids,
produces a presentation with parabolic generators
$\pmatrix{1& -1\cr 0& 1\cr}$
and
$\pmatrix{1& 0\cr 1& 1\cr}$
and describes the action of the derived group
on Poincar\'e's half plane.
}

\vskip10mm

Le groupe
 $SL(2, {\Bbb Z})$
 agit sur l'espace vectoriel
 $V\!=\!{\Bbb R}^{2}$%
 , respectant le r\'eseau
 $\Lambda\!=\!{\Bbb Z}^{2}$%
 , il induit donc un automorphisme du tore quotient
$\displaystyle {\Bbb T}\!=\!V/\Lambda$%
.
La suite exacte

$$0\rightarrow \Lambda\rightarrow V\rightarrow {\Bbb T}\rightarrow 0$$
fait appara\^{\i}tre
$\Lambda$
comme groupe fondamental de
${\Bbb T}$
et 
$\overline{M}: ({\Bbb T}, 0)\rightarrow ({\Bbb T}, 0)$%
, l'automorphisme induit par
$M\in SL(2, {\Bbb Z})$
produit sur ce groupe fondamental la restriction
$\pi_1(\overline{M}, 0) = M_{|} : \Lambda\rightarrow\Lambda$
de
$M$
\`a
$\Lambda$%
.

Le centre
$\{\pm I\} < SL(2, {\Bbb Z})$
fixe le sous-groupe
de
$2$%
-torsion
${\Bbb T}_2 = {{1}\over{2}}{\Bbb Z}/{\Bbb Z}$
et l'action de
$SL(2, {\Bbb Z})$
passe au quotient en une action par hom\'eomorphismes affines

$$\pi : G=PSL(2, {\Bbb Z}) \rightarrow
Affeo({\Bbb T}/\{\pm I\}, {\Bbb T}_2,  \{\overline{0}\})$$
du groupe modulaire
$G\!=\!SL(2, {\Bbb Z})/\{\pm I\}$
sur la sph\`ere plate \`a
$3 + 1$
 points singuliers (d'angle
$\pi$%
)%
\note{\eightpoint
 not\'es
$m = ({{1}\over{2}}, 0), n = (0, {{1}\over{2}}), p = m + n
= ({{1}\over{2}}, {{1}\over{2}})$
et
$\infty = (0, 0)$%
.
}
$({\Bbb T}/\{\pm I\}, {\Bbb T}_2, \{0\})
=({\Bbb S}, \{m, n, p, \infty\}, \{\infty\})=
({\Bbb S}, R, \{\infty\})$%
.

D'o\`u par lissage%
\note{\eightpoint
Si
$M=\pmatrix{a& b\cr c& d\cr}$
on d\'eforme
$\overline{M}$
au voisinage de
$R$
jusqu'\`a
$\overline{M}_c$%
, affine dans un plus petit voisinage
de partie lin\'eaire
$\pmatrix{a& -c\cr c& a\cr}$%
, la \og conformis\'ee sur la premi\`ere colonne\fg
de
$M$%
.
}%
, fibration de Cerf
et suite exacte  d'une fibration~:

$$\pi_0 : SL(2, {\Bbb Z})\rightarrow
\pi_0(Affeo({\Bbb S}, R,  \{\infty\}))\rightarrow$$
$$\rightarrow\pi_0(Diff({\Bbb S}, R,  \{\infty\}))
= \pi_1(Pl(R\setminus\{\infty\}, {\Bbb S}\setminus\{\infty\}))/Z$$
o\`u
$Z = <ababab>$
est le centre du groupe
$B_3$
des tresses \`a trois brins d'Artin%
\note{\eightpoint
$a$
(resp.
$b$%
)
 fixe
$n$
(resp.
$m$%
)
et  \'echange, par demi tour positif,
$m$
(resp.
$n$%
)
et
$p$%
.
}

$$B_3 = \pi_1(Pl(\{m, n, p,\}, {\Bbb S}\setminus\{\infty\}) =\
<a, b\ |\ aba = bab>$$

Si
$A= \pmatrix{1& -1\cr 0& 1\cr}, B = \pmatrix{1& 0\cr 1& 1\cr}$
et
$X = \pmatrix{0& -1\cr 1& 0\cr}$%
, on a~:

$$\pi_0(A) = a, \pi_0(B) = b\quad
\hbox{\rm et\/}\quad
ABA = X = BAB, X^{4} = I, X^{2} = -I\ne I$$
d'o\`u un morphisme
$\sigma : B_3\rightarrow SL(2, {\Bbb Z})$
d\'efini sur les g\'en\'erateurs de la pr\'esentation d'Artin par
$\sigma(a) = A, \sigma(b) = B$
qui, si on note
$\rho : B_3\rightarrow B_3/Z$%
, v\'erifie
$\pi_0\circ \sigma = \rho$%
.
\vskip2mm

Ainsi, si
$x = aba=bab$%
, le morphisme
$\sigma$
induit un morphisme injectif~:
$$\displaystyle\overline{\sigma} : B_3/\!\!<x^{4}>\ \rightarrow SL(2, {\Bbb Z})$$
qui est  isomorphisme puisqu'aussi surjectif car,  si 
$M\!\in\!SL(2, {\Bbb Z})$
et
$\mu\in B_3$
repr\'esente
$\pi_0(M)$%
, les matrices 
$M$
et
$\sigma(\mu)$
sont la
$\pi_1$%
-action de relev\'es de deux diff\'eomorphismes
homotopes de
$({\Bbb S}, R, \{\infty\}) = ({\Bbb T}/\{\pm I\}, {\Bbb T}_2, \{0\})$
donc sont, soit \'egales
et
$M = \sigma(\mu)$%
, soit oppos\'ees et
$M = \sigma(x^{2}\mu)$%
.\hfill\carre
\vskip2mm

On vient d'\'etablir les pr\'esentations~:
$$SL(2, {\Bbb Z}) = <A, B\ |\ ABA = BAB, (ABABAB)^{2}>$$
$$PSL(2, {\Bbb Z}) = <A, B\ |\ ABA = BAB, ABABAB>$$
et
$SL(2, {\Bbb Z})$
et
$PSL(2, {\Bbb Z})$
ont leurs ab\'elianis\'es cycliques d'ordre 
$12$
et $6$%
.\hfill\findem

\vfill\eject
Soit dans le demi-plan de Poincar\'e
${\Bbb H}=\{z\in{\Bbb C} | \Im(z)>0\}$
les translat\'es
par les puissances de
$A$
du domaine fondamental usuel de l'action de
$PSL(2, {\Bbb Z})$%
~:
$${\cal B}_{n} = A^{-n}\bigl(\{z\in{\Bbb C}\ |\
\Im(z)>0, -{{1}\over{2}}\leq \Re(z)\leq {{1}\over{2}}, |z|\geq 1\}\bigr)$$

L'union
$\tilde{H}\!=\!\cup {\cal B}_n$
quotient\'ee par l'image dans
$G$
du sous-groupe
$<\!\!A^{-6}X^{2}\!\!>$
du groupe d\'eriv\'e%
\note{\eightpoint
car,  dans l'ab\'enialis\'e
$\overline{A}=\overline{B}$
et
$\overline{X}=\overline{A}^{3}$%
, donc
$\overline{A^{-6}X^{2}}=\overline{1}$%
.}
$SL(2, {\Bbb Z})'$
 de
$SL(2, {\Bbb Z})$
est un hexagone hyperbolique point\'e qui,
par identification de ses c\^ot\'es oppos\'es,
produit un tore sym\'etrique point\'e
${\cal T}$%
, identifications par les images dans
$G$
des \'el\'ements suivants du groupe d\'eriv\'e~:

$$f_n = A^{-n}A^{-3}XA^{n} = A^{-(n + 3)}XA^{n}
= A^{-(n + 1)}A^{-1}BA^{(n + 1)}$$
$$f_{n + 1}f_{n - 1} = A^{-(n + 4)}XA^{(n + 1)}\!A^{-(n + 2)}\!XA^{n-1}
= A^{-(n + 3)}BAA^{(n + 1)}\!A^{-(n + 2)}XA^{n-1}$$
$$= A^{-(n + 3)}BXA^{n-1} = A^{-(n + 3)}XAA^{n-1} = A^{-(n + 3)}XA^{n} = f_n$$
De plus 
$A^{-6}X^{2}\!=\!f_{n} f_{n-3}\!=\!(f_{n-1}f_{-(n-2)}^{-1})(f_{n-1}^{-1}f_{n-2})\!=\!%
f_{n-1}[f_{n-2}, f_{n-1}]f_{n-1}^{-1}$%
.

Ainsi%
\note{\eightpoint
notant par abus
$f_{n}\in G$
l'images dans le groupe modulaire
$G$
de
$f_{n}\in SL(2, {\Bbb Z})$%
.}
${\cal T}$
est quotient de
${\Bbb H}$
par le sous-groupe de
$G$
engendr\'e par les
$f_{n}$%

$$<f_n; n\in{\Bbb Z}>\, =\, <f_{-2}, f_{-1}>\, =\, <f_{-1}^{-1}, f_{-2}>$$
qui est d'indice
$6$
dans le groupe modulaire
$G$
(l'indice de son l'ab\'elianis\'e
$G'$%
).

Comme
$<f_{n}; n\in{\Bbb Z}>\,\subset G'$%
, ce sous-groupe est le groupe d\'eriv\'e de
$G$%
~:
$$G'=PSL(2, {\Bbb Z})'=<f_{-2}, f_{-1}>$$
${\cal T}$
est le tore modulaire, de cusp de parabolique
le commutateur de
$f_{-1}^{-1}$ et
$f_{-2}$%
.\hfill\carre

Comme
$-I=X^{2}\in SL(2, {\Bbb Z})$
est d'ab\'elianis\'e
$\overline{A}^{6}\ne0,\ -I\not\in SL(2, {\Bbb Z})'$%
, ainsi le quotient
$SL(2, {\Bbb Z})\rightarrow G$
est injectif sur le groupe d\'eriv\'e
$SL(2, {\Bbb Z})'$ 
et ce groupe d\'eriv\'e est engendr\'e%
\note{\eightpoint librement puisque c'est le groupe fondamental d'un tore point\'e.
}
par
$f_{-2}=BA^{-1}=\pmatrix{1& 1\cr 1& 2\cr}$
et
$f_{-1}=A^{-1}B=\pmatrix{2& 1\cr 1& 1\cr}$%
.
\hfill\findem
\vfill\eject
\null
\vskip5mm
\centerline{\petcap Appendice}
\vskip7mm
\centerline{Description de
${\Bbb T}/\{\pm I\}$
et v\'erification de
$\pi_0(A)=a, \pi_0(B)=b$%
.
}
\vskip10mm
En identifiant
$(V, \Lambda)$
au plan complexe muni du r\'eseau de Gau\stz
$({\Bbb C}, {\goth G})$
o\`u
${\goth G}={\Bbb Z}+i{\Bbb Z}$%
, la fonction
$\wp : {\Bbb C}\rightarrow{\Bbb C}$
de Weierstra\stz~:
$$\wp(z)={{1}\over{z^{2}}}+\sum_{\gamma\in\Gamma\setminus\{0\}}
{{1}\over{(z-\gamma)^{2}}} -{{1}\over{\gamma^{2}}}$$
qui est paire et invariante par
$\Gamma$%
, identifie
${\Bbb T}/\{\pm I\}$
\`a
${\Bbb C}\cup\{\infty\}$
avec~:
$$\wp(0)=\infty,\
\wp({{1}\over{2}}(1+i))=0,\
\wp({{1}\over{2}})=-\wp({{1}\over{2}}i)\in\, ]0, +\infty[$$

Les deux \og sous-groupes de coordonn\'es\fg
${\bf t}={\Bbb R}/{\Bbb Z},\ {\bf u}={\Bbb R}i/{\Bbb Z}i$
vont%
\note{\eightpoint
pour le v\'erifier il suffit de remarquer que, pr\`es de
$0$%
, on a
$\wp(z)\sim{{1}\over{z^{2}}}$%
.}
sur les demi-axes r\'eels tronqu\'es
$\wp({\bf t})=H_{\wp(m)}^{+}:=[\wp(m), +\infty[\cup\{\infty\}\subset{\Bbb R}\cup\{\infty\}$
et
$\wp({\bf u})=H^{\wp(n)}_{-}:=]-\infty, \wp(n)]\cup\{\infty\}\subset{\Bbb R}\cup\{\infty\}$
et les deux
\og sous-groupes diagonaux\fg
$\Delta_{+}={\Bbb R}(1+i)/{\Bbb Z}(1+i)$
et
$\Delta_{-}={\Bbb R}(1-i)/{\Bbb Z}(1-i)$
vont sur les demi-axes imaginaires
$\wp(\Delta_{+})=V^{+}:=[0, +\infty[\,i\cup \{\infty\}$ 
et
$\wp(\Delta_{-})=V_{-}:=]-\infty, 0]\,i\cup\{\infty\}$
.

Ainsi l'action de
$A=\pmatrix{1& -1\cr 0& 1\cr}$
fixe
$H^{+}_{\wp(m)}$
et envoie
$H_{-}^{\wp(n)}$
et
$V^{+}$
sur
$V_{-}$
et
$H_{-}^{\wp(n)}$
respectivement, c'est le \og demi-tour\fg positif
$a$%
.

Et l'action de
$B=\pmatrix{1& 0\cr 1& 1\cr}$
fixe
$H_{-}^{\wp(n)}$
et envoie
$H^{+}_{\wp(m)}$
et
$V_{-}$
sur
$V^{+}$
et
$H_{-}^{\wp(n)}$
respectivement, c'est le \og demi-tour\fg positif
$b$%
.
\hfill\findem

\vskip15mm
\centerline{{\petcap Remerciements}}

Cette note doit tout aux relectures attentives de Danielle Bozonat, Daniel Marin
et surtout Greg McShane qui,
en partageant ses tentatives contre
la conjecture d'unicit\'e de Frobenius pour les nombres de Markoff
(Cf. {\bf [H]}),
a pouss\'e le r\'edacteur \`a comprendre le tore modulaire,
sans oublier l'objectivit\'e de l'arbitre dont le rapport a fourni la mati\`ere
des commentaires bibliographiques ci-dessous.

\vfill\eject
\centerline{{\petcap Commentaires bibliographiques}}
\vskip3mm

Les seuls calculs
$ABA\!=\!BAB, f_{n+1}f_{n-1}\!=\!f_{n}$
et
$\pi_{0}(A)\!=\!a, \pi_{0}(B)\!=\!b$
plut\^ot que ceux et les  r\'ecurrences \`a base d'algorithme d'Euclide
du traitement classique (p. 9-19 de {\bf [R]})
profitent de ce que  la suffisance des relations est donn\'ee par
la pr\'esentation du groupe des 
classes d'isotopie de diff\'eomorphismes de la sph\`ere plate
$({\Bbb S}, \{m, n, p, \infty\}, \{\infty\})=
({\Bbb T}/\{\pm I\}, {\Bbb T}_{2}, \{0\})$
respectant points d'ordre
$2$
et
$1$%
, pr\'esentation d\'eduite de celle du
groupe des tresses \`a trois brins  d'Artin ({\bf [A]}).

Pr\`es d'un quart de si\`ecle plus tard,%
\note{\eightpoint
{\sl Theory of Braids} Ann. of Math. 48 (1947), p. 101-126.
}
Artin conc\'edait que (dans son article de 1925) ``{\sl Most of the proofs are entirely intuitive}'',
renvoyant pour les relations de
$B_{n}$
au quasiment simultan\'e%
\note{\eightpoint
{\sl The algebraic braid group.} Ann. Math. 48 (1947), p. 127-136.
}
traitement alg\'ebrique de F. Bohnenblust via l'action de
$B_{n}$
sur le groupe libre
$F_{n}$
(le groupe fondamental du plan priv\'e de
$n$
points).

Depuis l'\'etablissement de ces relations%
\note{\eightpoint
tant celles du groupe des tresses que le passage au groupe des classes d'isotopies
d'hom\'eo\-mor\-phismes ou diff\'eomorphismes respectant une partie finie
de la sph\`ere.
}
a fait couler beaucoup d'encre.

Un traitement clair et rigoureux est donn\'e par l'\'ecole d'Orsay
via la m\'ethode des espaces fonctionnels%
\note{\eightpoint
l'espace des plongements de
$n$
points dans le plan est stratifi\'e avec une seule strate ouverte qui est contractile,
les g\'en\'erateurs du groupe fondamental correspondant aux travers\'ees des states de codimension
$1$
et  les  relations aux contournements des strates de codimension
$2$%
.}
et les th\'eor\`emes de fibration de Cerf ({\bf [C]}).%
\note{\eightpoint
Cf. V. Po\'enaru Expos\'e 2, p. 21-31 de {\bf [FLP]}.
}

\vskip8mm
\centerline{{\petcap R\'ef\'erences}}

\vskip2mm
{\hangindent=1cm\hangafter=1\noindent{\bf   [A]}\
{\petcap Artin E}\pointir
 {\sl Theorie der Z\H opfe},
Hamb. Abh. 4,  p. 47-72%
, 
{\oldstyle 1925}.
\par}


\vskip1mm
{\hangindent=1cm\hangafter=1\noindent{\bf   [C]}\
{\petcap Cerf J}\pointir
 {\sl Topologie de certains espaces de plongements},\hfill\break
Bull. Soc. Math. France 89, p. 227-280 (Th\`ese Sc. math. Paris 1960)%
, 
{\oldstyle 1961}.
\par}

\vskip1mm
{\hangindent=1cm\hangafter=1\noindent{\bf   [FLP]}\
{\petcap Fathi A., Laudenbach F., Po\'enaru V}\pointir
 {\sl Travaux de Thurston sur les surfaces.  S\'eminaire Orsay} 
Ast\'erisque 66-67%
, 
{\oldstyle 1979}.
\par}

\vskip1mm
{\hangindent=1cm\hangafter=1\noindent{\bf   [H]}\
{\petcap Haas A}\pointir
 {\sl The geometry of Markoff forms}, dans
 {\sl Number theory} \hfill\break
Lecture Notes in Mathematics 1240, p. 135-144%
, 
{\oldstyle 1985}.
\par}

\vskip1mm
{\hangindent=1cm\hangafter=1\noindent{\bf   [R]}\
{\petcap Rankin R}\pointir
 {\sl Modular forms and functions},
C.U.P.,  p. 9-19%
, 
{\oldstyle 1977}.
\par}

\vskip10mm

\noindent{
Secr\'etaire et r\'edacteur~: Alexis Marin Bozonat\hfill\break
courriel~: alexis.charles.marin@gmail.com}

\end

\vskip3mm
{\hangindent=1cm\hangafter=1\noindent{\bf   [C2]}\
{\petcap Cerf J}\pointir
 {\sl Sur les diff\'eomorphismes de la sph\`ere de dimension trois (%
 $\Gamma_4=0$%
 )},\hfill\break
Lecture Notes in Math. 53, Springer%
, 
{\oldstyle 1968}.
\par}

\vskip3mm
{\hangindent=1cm\hangafter=1\noindent{\bf   [C2]}\
{\petcap Cerf J}\pointir
 {\sl La stratification naturelle des espaces de fonctions diff\'erentiables r\'eelles et le th\'eor\`eme de la pseudo-isotopie},\hfill\break
Pub. math. I.H.E.S. 39, p.5-173%
, 
{\oldstyle 1970}.
\par}

\vskip3mm
{\hangindent=1cm\hangafter=1\noindent{\bf   [S]}\
{\petcap Smale S}\pointir
 {\sl Diffeomorphisms of the 2-sphere},\hfill\break
Proc. Amer. Math. Soc. 10, p. 621-626%
, 
{\oldstyle 1959}.
\par}

\end

L'action de
$\{\pm I\}$
sur
${\Bbb T}$
a pour  domaine fondamental le rectangle
$[0, 1]\!\times\![0, {{1}\over{2}}]$%
.

Les c\^ot\'es
$0n\!=\!\{0\}\!\times\![0, {{1}\over{2}}]$
et
$0'n'\!=\!\{1\}\!\times\![0, {{1}\over{2}}]$
sont identifi\'es par translation.

Les identifications sur les bases inf\'erieures
$00'=[0, 1]\times\{0\}$
et sup\'erieures
$nn'=[0, 1]\times\{{{1}\over{2}}\}$
sont les pliages en leurs milieux
$m$
et
$p$
respectivement.

En prenant pour mod\`ele aux points d'ordre
$2$
de l'application quotient
$(x, y)\mapsto (x^2 - y^2, 2xy)$
(le carr\'e complexe) qui double les angles au centre et pensant le point
$0$
\`a l'infini, on repr\'esente
${\Bbb T}/\{\pm I\}\setminus\{0\}$
comme le plan
${\Bbb R}^2$%
.

Le point
$p$
s'identifie au centre.\hfill\break
Les segments
$0m$
 au \og d\'ebut d'axe vertical\fgf~:
$V^{m}_{-}=\{0\}\times]-\infty, -1]$%
.\hfill\break
Les segments
$0n$
et
$0n'$
 \`a la \og fin d'axe vertical\fgf~:
$V_{n}^{+}=\{0\}\times[1, +\infty[$%
.\hfill\break
Les deux diagonales
$0p$
et
$p0'$
(correspondant aux sous-groupes
${\Bbb R}(1, 1)/{\Bbb Z}(1, 1)$
et
${\Bbb R}(1, -1)/{\Bbb Z}(1, -1)$
de
${\Bbb T}={\Bbb R}^{2}/{\Bbb Z}^{2}$
s'identifient aux demi-axes horizontaux n\'egatifs et positifs
$H_{-}^{0}=]\!-\infty, 0]\times\{0\}$
et
$H_{0}^{+}=[0, +\infty[\times\{0\}$
respectivement.

Ainsi l'action de
$A=\pmatrix{1& 0\cr -1& 1\cr}$
fixe
$V_{n}^{+}$
et envoie
$V_{-}^{m}$
sur
$H_{0}^{+}$%
, c'est
$a$%
.

L'action de
$B=\pmatrix{1& 1\cr 0& 1\cr}$
fixe
$V_{-}^{m}$
et envoie
$V_{n}^{+}$
sur
$H_{-}^{0}$%
, c'est
$b$%
.

%% file: JOSEPH.tex
\input option_keys
\let\stz\ss
\input SSIND
\input FONTESMATH

\input MATH

\input pdffix.tex

%% file: option_keys
\ifx\optionkeymacros\undefined\else \fi

\catcode`\Œ=\active\defŒ{{\aa}}       
\catcode`\º=\active\defº{\int}        
\catcode`\=\active\def{\c c}        
\catcode`\¶=\active\def¶{\partial}    
\catcode`\Ä=\active\defÄ{\oint}       
\catcode`\Æ=\active\defÆ{\triangle}   
\catcode`\Â=\active\defÂ{\neg}        
\catcode`\µ=\active\defµ{\mu}         
\catcode`\¿=\active\def¿{{\o}}        
\catcode`\¹=\active\def¹{\pi}         
\catcode`\Ï=\active\defÏ{{\oe}}       
\catcode`\§=\active\def§{{\ss}}       
\catcode`\ =\active\def {\dagger}     
\catcode`\Ã=\active\defÃ{\sqrt}       
\catcode`\·=\active\def·{\Sigma}      
\catcode`\Å=\active\defÅ{\approx}     
\catcode`\½=\active\def½{\Omega}      
\catcode`\£=\active\def£{{\it\$}}     
\catcode`\°=\active\def°{\infty}      
\catcode`\¤=\active\def¤{{\S}}        
\catcode`\¦=\active\def¦{{\P}}        
\catcode`\¥=\active\def¥{\bullet}     
\catcode`\»=\active\def»{\leavevmode\raise.585ex\hbox{\b a}}      
\catcode`\¼=\active\def¼{\leavevmode\raise.6ex\hbox{\b o}}        
\catcode`\­=\active\def­{\not=}       
\catcode`\²=\active\def²{\leq}        
\catcode`\³=\active\def³{\geq}        
\catcode`\Ö=\active\defÖ{\div}        
\catcode`\É=\active\defÉ{{\dots}}     
\catcode`\¾=\active\def¾{{\ae}}       
\catcode`\Ç=\active\defÇ{\ll}         
\catcode`\Ò=\active\defÒ{``}          
\catcode`\Á=\active\defÁ{!`}          
\catcode`\¢=\active\def¢{\rlap/c}     


\catcode`\=\active\def{{\AA}}       
\catcode`\'=\active\def'{\c C}        
\catcode`\¯=\active\def¯{{\O}}        
\catcode`\¸=\active\def¸{\Pi}         
\catcode`\Î=\active\defÎ{{\OE}}       
\catcode`\®=\active\def®{{\AE}}       
\catcode`\×=\active\def×{\diamond}    
\catcode`\¡=\active\def¡{\accent'27}  
\catcode`\Ó=\active\defÓ{''}          
\catcode`\±=\active\def±{\pm}         
\catcode`\È=\active\defÈ{\gg}         
\catcode`\À=\active\defÀ{?`}          
\catcode`\Ð=\active\defÐ{--}          
\catcode`\Ñ=\active\defÑ{---}         


\catcode`\Š=\active\defŠ{\"a}        
\catcode`\'=\active\def'{\"e}        
\catcode`\•=\active\def•{\"{\i}}     
\catcode`\š=\active\defš{\"o}        
\catcode`\Ÿ=\active\defŸ{\"u}        
\catcode`\Ø=\active\defØ{\"y}        
\catcode`\€=\active\def€{\"A}        
\catcode`\…=\active\def…{\"O}        
\catcode`\†=\active\def†{\"U}        
\catcode`\‡=\active\def‡{\'a}        
\catcode`\Ž=\active\defŽ{\'e}        
\catcode`\'=\active\def'{\'{\i}}     
\catcode`\—=\active\def—{\'o}        
\catcode`\œ=\active\defœ{\'u}        
\catcode`\ƒ=\active\defƒ{\'E}        
\catcode`\ˆ=\active\defˆ{\`a}        
\catcode`\=\active\def{\`e}        
\catcode`\"=\active\def"{\`{\i}}     
\catcode`\˜=\active\def˜{\`o}        
\catcode`\=\active\def{\`u}        
\catcode`\Ë=\active\defË{\`A}        
\catcode`\‹=\active\def‹{\~a}        
\catcode`\–=\active\def–{\~n}        
\catcode`\›=\active\def›{\~o}        
\catcode`\Ì=\active\defÌ{\~A}        
\catcode`\"=\active\def"{\~N}        
\catcode`\Í=\active\defÍ{\~O}        
\catcode`\‰=\active\def‰{\^a}        
\catcode`\=\active\def{\^e}        
\catcode`\"=\active\def"{\^{\i}}     
\catcode`\™=\active\def™{\^o}        
\catcode`\ž=\active\defž{\^u}        

\let\optionkeymacros\null

%% file: SSIND.tex
\font\eightrm=cmr8
\font\eighti=cmmi8
\font\eightsy=cmsy8
\font\eightbf=cmbx8
\font\eighttt=cmtt8
\font\eightit=cmti8
\font\eightsl=cmsl8
\font\sixrm=cmr6
\font\sixi=cmmi6
\font\sixsy=cmsy6
\font\sixbf=cmbx6
\skewchar\eighti='177 \skewchar\sixi='177
\skewchar\eightsy='60 \skewchar\sixsy='60
\catcode`\@=11

\def\pc#1#2|{{\bigf@ntpc #1\penalty\@MM\hskip\z@skip\smallf@ntpc #2}}
\def\tenpoint{%
\textfont0=\tenrm \scriptfont0=\sevenrm \scriptscriptfont0=\fiverm
\def\rm{\fam\z@\tenrm}%
\textfont1=\teni \scriptfont1=\seveni \scriptscriptfont1=\fivei
\def\oldstyle{\fam\@ne\teni}%
  \textfont2=\tensy \scriptfont2=\sevensy \scriptscriptfont2=\fivesy
  \textfont\itfam=\tenit
  \def\it{\fam\itfam\tenit}%
  \textfont\slfam=\tensl
  \def\sl{\fam\slfam\tensl}%
  \textfont\bffam=\tenbf \scriptfont\bffam=\sevenbf
  \scriptscriptfont\bffam=\fivebf
  \def\bf{\fam\bffam\tenbf}%
  \textfont\ttfam=\tentt
  \def\tt{\fam\ttfam\tentt}%
  \abovedisplayskip=6pt plus 2pt minus 6pt
  \abovedisplayshortskip=0pt plus 3pt
  \belowdisplayskip=6pt plus 2pt minus 6pt
  \belowdisplayshortskip=7pt plus 3pt minus 4pt
  \smallskipamount=3pt plus 1pt minus 1pt
  \medskipamount=6pt plus 2pt minus 2pt
  \bigskipamount=12pt plus 4pt minus 4pt
  \normalbaselineskip=12pt
  \setbox\strutbox=\hbox{\vrule height8.5pt depth3.5pt width0pt}%
  \let\bigf@ntpc=\tenrm \let\smallf@ntpc=\sevenrm
  \normalbaselines\rm}
\def\eightpoint{%
  \textfont0=\eightrm \scriptfont0=\sixrm \scriptscriptfont0=\fiverm
  \def\rm{\fam\z@\eightrm}%
  \textfont1=\eighti \scriptfont1=\sixi \scriptscriptfont1=\fivei
  \def\oldstyle{\fam\@ne\eighti}%
  \textfont2=\eightsy \scriptfont2=\sixsy \scriptscriptfont2=\fivesy
  \textfont\itfam=\eightit
  \def\it{\fam\itfam\eightit}%
  \textfont\slfam=\eightsl
  \def\sl{\fam\slfam\eightsl}%
  \textfont\bffam=\eightbf \scriptfont\bffam=\sixbf
  \scriptscriptfont\bffam=\fivebf
  \def\bf{\fam\bffam\eightbf}%
  \textfont\ttfam=\eighttt
  \def\tt{\fam\ttfam\eighttt}%
  \abovedisplayskip=9pt plus 2pt minus 6pt
  \abovedisplayshortskip=0pt plus 2pt
  \belowdisplayskip=9pt plus 2pt minus 6pt
  \belowdisplayshortskip=5pt plus 2pt minus 3pt
  \smallskipamount=2pt plus 1pt minus 1pt
  \medskipamount=4pt plus 2pt minus 1pt
  \bigskipamount=9pt plus 3pt minus 3pt
  \normalbaselineskip=9pt
  \setbox\strutbox=\hbox{\vrule height7pt depth2pt width0pt}%
  \let\bigf@ntpc=\eightrm \let\smallf@ntpc=\sixrm
  \normalbaselines\rm}

\newskip\LastSkip
\def\nobreakatskip{\relax\ifhmode\ifdim\lastskip>\z@
  \LastSkip\lastskip\unskip\nobreak\hskip\LastSkip
  \fi\fi}
\catcode`\;=\active
\catcode`\:=\active
\catcode`\!=\active
\catcode`\?=\active
\def;{\nobreakatskip\string;}
\def:{\nobreakatskip\string:}
\def!{\nobreakatskip\string!}
\def?{\nobreakatskip\string?}
%
  \let\titlefont=\seventeenrm

\font\tenss=cmss10 \font\tencaps=cmcsc10
\newfam\ssfam   \textfont\ssfam=\tenss     \def\ss{\fam\ssfam\tenss}
\newfam\capsfam \textfont\capsfam=\tencaps \def\petcap{\fam\capsfam\tencaps}

\font\teneuf=eufm10 \font\seveneuf=eufm7 \font\fiveeuf=eufm5
\font\tenmsa=msam10 \font\sevenmsa=msam7 \font\fivemsa=msam5
\font\tenmsb=msbm10 \font\sevenmsb=msbm7 \font\fivemsb=msbm5

\newfam\msafam \textfont\msafam=\tenmsa \scriptfont\msafam=\sevenmsa
\scriptscriptfont\msafam=\fivemsa
\newfam\msbfam \textfont\msbfam=\tenmsb \scriptfont\msbfam=\sevenmsb
\scriptscriptfont\msbfam=\fivemsb \def\Bbb{\fam\msbfam }
\newfam\euffam \textfont\euffam=\teneuf \scriptfont\euffam=\seveneuf
\scriptscriptfont\euffam=\fiveeuf \def\goth{\fam\euffam\teneuf}

\newskip\LastSkip
\def\nobreakatskip{\relax\ifhmode\ifdim\lastskip>0pt
  \LastSkip\lastskip\unskip
  \nobreak\hskip\LastSkip
  \fi\fi}
\catcode`\;=\active \def;{\nobreakatskip\string;}
\catcode`\:=\active \def:{\nobreakatskip\string:}
\catcode`\!=\active \def!{\nobreakatskip\string!}
\catcode`\?=\active \def?{\nobreakatskip\string?}

\frenchspacing
\tenpoint


\magnification 1200

\pretolerance=500   \tolerance=1000   
\brokenpenalty=5000 

\hoffset=1cm
\hsize=125mm \vsize=187mm \parindent=1cm

\parskip=5pt plus 1pt
\newif\ifpagetitre 				\pagetitretrue
\newtoks\hautpagetitre			\hautpagetitre={\hfil}
\newtoks\baspagetitre       \baspagetitre={\hfil\tenrm\folio\hfil}
\newtoks\auteurcourant     \auteurcourant={\hfil}
\newtoks\titrecourant       \titrecourant={\hfil}
\newtoks\hautpagegauche  \newtoks\hautpagedroite
\newtoks\chapitre              \newtoks\paragraphe
\hautpagegauche={{\tenrm\folio}\hfil\the\auteurcourant\hfil}
\hautpagedroite={\hfil\the\titrecourant\hfil{\tenrm\folio}}
\newtoks\baspagegauche     \baspagegauche={\hfil}
\newtoks\baspagedroite     \baspagedroite={\hfil}
\headline={\ifpagetitre\the\hautpagetitre
\else\ifodd\pageno\the\hautpagedroite
\else\the\hautpagegauche\fi\fi}
\footline={\ifpagetitre\the\baspagetitre
\global\pagetitrefalse
\else\ifodd\pageno\the\baspagedroite
\else\the\baspagegauche\fi\fi}

\catcode`\@=11
\def\p@int{{\rm .}}
\def\p@intir{\discretionary{\rm .}{}{\rm .\kern.35em---\kern.7em}}
\def\pointir{\afterassignment\pointir@\let\next=}
\def\pointir@{\ifx\next\par\p@int\else\p@intir\fi\next}
\catcode`\@=12

\def\carre{\hbox{\tensy\char'164\hskip-6.66666pt\char'165}}

\def\cite#1{[#1]}
\def\buildo#1\over#2{\mathrel{\mathop{\null#2}\limits^{#1}}}
\def\buildu#1\under#2{\mathrel{\mathop{\null#2}\limits_{#1}}}

\def\_#1{_{\vtop{\halign{$\scriptstyle{##}$&
                       $\scriptstyle{{}##}$\cr #1\crcr}}}}

\def\text#1{\hbox{\rm #1}} 

\def\article#1|#2|#3|#4|#5|#6|#7|
  {\leftskip=7mm\noindent
   \hangindent=2mm\hangafter=1
    \llap{\bf [#1]\hskip1.35em}{\petcap #2}\pointir {\sl #3}, {\rm #4},
    \nobreak {\bf #5} ({\oldstyle #6}), \nobreak #7.\par}
\def\livre#1|#2|#3|#4|#5|
 {\leftskip=7mm\noindent
   \hangindent=2mm\hangafter=1
    \llap{\bf [#1]\hskip1.35em}{\petcap #2}\pointir {\sl #3}, {\rm #4},
    {\oldstyle #5}.\par}
\def\alivre#1|#2|#3|#4|#5|#6|#7|
  {\leftskip=7mm\noindent
    \hangindent=2mm\hangafter=1
    \llap{\bf [#1]\hskip1.35em}{\petcap #2}\pointir {\sl #3},
dans {\it #4}, \'edit\'e par {\bf #5}, {\rm #6}, {\oldstyle #7}.\par}
\def\divers#1|#2|#3|#4|
  {\leftskip=7mm\noindent
   \hangindent=2mm\hangafter=1
    \llap{[#1]\hskip1.35em}{\petcap #2}\pointir #3,
    {\oldstyle #4}.\par}
\newcount\secno \secno=0
\newcount\ssecno \ssecno=0
\newcount\sssecno \sssecno=0
\newcount\chapno \chapno=0
\newcount\notenumber \notenumber=1
\newcount\exino \exino=0
\newcount\expno \expno=0
%
%
%
\newdimen\indentsec\indentsec=20pt
\newdimen\indentssec\indentssec=20pt
\newdimen\indentsssec\indentsssec=20pt
\newdimen\indentrem\indentssec=0pt

%
\newdimen\indentTh\indentTh=0pt

%
\newdimen\indentth\indentssec=0pt
\def\titregen#1#2|{\par\vskip .5cm\penalty -100
             {\leftskip=0pt plus \hsize
			   \rightskip=\leftskip
			   \parfillskip=0pt
			   \baselineskip=17pt
			   \noindent
             \titlefont{ #1}{#2}\par}
			  \vskip 5pt\penalty 500}
\def\titre#1|{\titregen{}{#1}|}
%
\def\ntitre#1|{\global\advance\chapno by 1\global\ssecno=0\global\secno=0
\titregen{\the\chapno\ }{#1}|}
\def\auteur#1|{\vskip 5pt\penalty 100
               \vbox{\centerline{\sl par #1}
                    \vskip 5pt}\penalty 100}
\def\sectiongen#1#2#3{\parindent=\indentsec\par\vskip .3cm
\vskip 0mm plus -20mm minus 1,5mm\penalty-50
{\bf #1}{\bf #2}{#3}\nobreak\parindent=20pt}%
%
\def\secc#1|{\sectiongen{}{#1}{\pointir}}
%
\def\nsecc#1|%
{\global\advance\secno by 1\global\ssecno=0\global\sssecno=0
\sectiongen{\the\secno\ }{#1}{\pointir}}
%
%
\def\secp#1|{\sectiongen{}{#1}{}\par}
%
\def\nsecp#1|%
{\global\advance\secno by 1\global\ssecno=0\global\sssecno=0
\sectiongen{\the\secno\ }{#1}{}\par}
\def\ssectiongen#1#2#3#4{\parindent=\indentssec\par\vskip .2cm
\vskip 0mm plus -20mm minus 1,5mm\penalty-50
{\bf #1}{\sl #2}{\sl #3}{#4}\nobreak\medskip\parindent=20pt}%
%
\def\ssecc#1|#2{\ssectiongen{}{#1}{\pointir}{#2}}
%
\def\nssecc#1|#2{\global\advance\ssecno by 1\global\sssecno=0
\ssectiongen{\the\secno.\the\ssecno\ }{#1}{\pointir}{#2}}
%
\def\ssecp#1|{\ssectiongen{}{#1}{}{}\par}
%
\def\nssecp#1|{\global\advance\ssecno by 1\global\sssecno=0
\ssectiongen{\the\secno.\the\ssecno\ }{#1}{}{}\par}
\def\sssectiongen#1#2#3#4{\parindent=\indentsssec\par\vskip .2cm
\vskip 0mm plus -20mm minus 1,5mm\penalty-50
{\bf #1}{\sl #2}{\sl #3}{#4}\nobreak\medskip\parindent=20pt}%
%
\def\sssecc#1|#2{\sssectiongen{}{#1}{\pointir}{#2}}
%
\def\nsssecc#1|#2{\global\advance\sssecno by 1
\ssectiongen{\the\secno.\the\ssecno.\the\sssecno\ }{#1}{\pointir}{#2}}
%
\def\sssecp#1|{\sssectiongen{}{#1}{}{}\par}
%
\def\nsssecp#1|{\global\advance\sssecno by 1
\sssectiongen{\the\secno.\the\ssecno.\the\sssecno\ }{#1}{}{}\par}
\long\def\Th#1#2#3#4{\parindent=\indentTh\par\vskip5pt
{#1}{\petcap #2}{\sl #3}\parindent=20pt{\sl #4\par}\vskip 5pt\parindent=20pt}
\long\def\th#1#2#3#4{\parindent=\indentth\par\vskip5pt
{#1}{\pppetcap #2}{\eightpoint\sl #3}\parindent=20pt{\eightpoint\sl #4\par}\vskip 5pt\parindent=20pt}
\long\def\remarque#1#2#3#4{\parindent=\indentrem\par\vskip5pt
{#1}{\eightpoint \sl #2}{\eightpoint\sl #3}\parindent=20pt{\eightpoint#4\par}
\vskip 5pt\parindent=20pt}
\long\def\remarques#1#2#3#4{\parindent=\indentrem\par\vskip5pt
{#1}{\eightpoint \sl #2}{\eightpoint\sl #3}{\eightpoint#4\par}
\vskip 5pt\parindent=20pt}

\long\def\remarquesa#1#2#3#4#5{\parindent=\indentrem\par\vskip5pt
{#1}{\eightpoint \sl #2}{\eightpoint\sl #3}{\eightpoint#4}
\parindent=20pt{\eightpoint#5\par}
\vskip 5pt\parindent=20pt}

\long\def\Remarque#1#2#3#4{\parindent=\indentrem\par\vskip5pt
{#1}{ \sl #2}{\sl #3}\parindent=20pt{#4\par}
\vskip 5pt\parindent=20pt}
\long\def\remarquesn#1#2#3#4{\parindent=\indentrem\par\vskip5pt
{\eightpoint \sl #1}{\eightpoint #2}{\eightpoint\sl #3}{\eightpoint#4\par}
\vskip 5pt\parindent=20pt}

\long\def\Remarquen#1#2#3#4{\parindent=\indentrem\par\vskip5pt
{ \sl #1}{#2}{\sl #3}\parindent=20pt{#4\par}
\vskip 5pt\parindent=20pt}

%
\long\def\Thc#1|#2\finc{\Th{}{#1}{\pointir}{#2}}
%
\long\def\Thnc#1|#2|#3\finnc{\Th{#1}{#2}{\pointir}{#3}}
\long\def\Exic#1|#2\finc{{\global\advance\exino by 1}\Remarquen%
{Exercice }{\the\exino}{ #1\pointir}{#2}}
\long\def\Expc#1|#2\finc{{\global\advance\expno by 1}\Remarquen%
{Exemple }{\the\expno}{ #1\pointir}{#2}}
\long\def\exic#1|#2\finc{{\global\advance\exino by 1}\remarquesn%
{\bf Exercice }{\the\exino}{ #1\pointir}{#2}}
\long\def\expc#1|#2\finc{{\global\advance\expno by 1}\remarquesn%
{Exemple }{\the\expno}{ #1\pointir}{#2}}

%
\long\def\Ec#1\finc{\Th{}{}{}{#1}}
%

\long\def\thc#1|#2\finc{\th{}{#1}{\pointir}{#2}}%
\long\def\thc#1\finc{\th{}{}{}{#1}}

\long\def\Thp#1|#2\finp{\Th{}{#1}{\par}{#2}}
%
\long\def\thp#1|#2\finp{\th{}{#1}{\par}{#2}}
%
\long\def\rmc#1|#2\finc{\remarque{}{#1}{\pointir}{#2}}
%
\long\def\Rmc#1|#2\finc{\Remarque{}{#1}{\pointir}{#2}}

\long\def\rmp#1|#2\finp{\remarque{}{#1}{\par}{#2}}
%
\long\def\Rmp#1|#2\finp{\Remarque{}{#1}{\par}{#2}}

\long\def\parc#1\finc{\remarque{}{}{}{#1}}
%
\long\def\parcs#1\fincs{\remarques{}{}{}{#1}}

\long\def\parcsa#1\fins#2\fincsa{\remarquesa{}{}{}{#1}{#2}}

\def\Rm#1|{\parindent=0pt\par\vskip5pt{\sl #1}\pointir\parindent=20pt}
%

%
\def\Demd#1|{\parindent=0pt\par{\sl D\'emonstration d#1}\pointir\parindent=20pt}
%


%

%

\def\preuved#1|{\parindent=0pt\par{\sl Preuve d#1}\pointir\parindent=20pt}


\def\qed{\quad\hbox{\hskip 1pt\vrule width 4pt height 6pt
         depth 1.5pt\hskip 1pt}}
\def\findem{\penalty 500 \hbox{\qed}\par\vskip 3pt}

%

%
%

%

%

%

%

%
\def\og{\leavevmode\raise.3ex\hbox{$\scriptscriptstyle\langle\!\langle\,$}}

\def\fgf{\/\leavevmode\raise.3ex\hbox{$\scriptscriptstyle\,\rangle\!\rangle$}}
%
\def\fg{\fgf\ }
\def\note#1{\footnote{$^{\the\notenumber}$}{#1}%
\global\advance\notenumber by 1}%
\def\ieme{\raise 1ex\hbox{\pc{}i\`eme|}\ }
\def\iemes{\raise 1ex\hbox{\pc{}i\`emes|}\ }

%% file: FONTESMATH.tex
\input amsfonts.def

%% file: MATH
\def\build#1_#2^#3{\mathrel{\mathop{\kern 0pt#1}\limits_{#2}^{#3}}}
%

\def\hdfl#1#2{\smash{\mathop{\hbox to 12mm{\rightarrowfill}}
\limits^{\scriptstyle#1}_{\scriptstyle#2}}}

\def\hdhfl#1#2{\smash{\mathop{\hbox to 12mm{\hookrightarrowfill}}
\limits^{\scriptstyle#1}_{\scriptstyle#2}}}

\def\hgfl#1#2{\smash{\mathop{\hbox to 12mm{\leftarrowfill}}
\limits^{\scriptstyle#1}_{\scriptstyle#2}}}

\def\hghfl#1#2{\smash{\mathop{\hbox to 12mm{\hookleftarrowfill}}
\limits^{\scriptstyle#1}_{\scriptstyle#2}}}













\def\th{\mathop{\rm th}\nolimits}




%% file: pdffix.tex
\input ifpdf.sty
\ifpdf
\pdfpagewidth=8.5truein
\pdfpageheight=11truein
\pdfhorigin=1truein
\pdfvorigin=1truein
\fi